\documentclass[12pt]{article}%
\usepackage{amsfonts}
\usepackage{amsmath}
\usepackage{amssymb}
\usepackage{color}
\usepackage{graphicx}%
\providecommand{\U}[1]{\protect\rule{.1in}{.1in}}
\newtheorem{theorem}{Theorem}

\newtheorem{corollary}[theorem]{Corollary}

\newtheorem{definition}[theorem]{Definition}

\newtheorem{lemma}[theorem]{Lemma}

\newtheorem{remark}[theorem]{Remark}

\newenvironment{proof}[1][Proof]{\noindent\textbf{#1.} }{\ \rule{0.5em}{0.5em}}
\begin{document}

\title{Monotone Multivalued Nonautonomous Dynamical Systems}
\author{Jos\'{e} A. Langa$^{1}$, Jacson Simsen$^{2}$, Mariza Stefanello Simsen$^{2}$,
Jos\'{e} Valero$^{3}$\\$^{1}${\small Dpto. Ecuaciones Diferenciales y An\'{a}lisis Num\'{e}rico}\\{\small Universidad de Sevilla, Apdo. de Correos 1160, 41080-Sevilla, Spain}\\{\small E.mail: langa@us.es}\\$^{2}${\small Instituto de Matem\'{a}tica e Computac\~{a}o, Universidade
Federal de Itajub\'{a},}\\\ \ {\small Av. BPS n. 1303, Bairro Pinheirinho, 37500-903, Itajub\'{a} - MG -
Brazil}\\\ \ {\small E.mails: jacson@unifei.edu.br;\ mariza@unifei.edu.br}\\$^{3}${\small Centro de Investigaci\'{o}n Operativa, Universidad Miguel
Hern\'{a}ndez de Elche, }\\{\small Avda. de la Universidad, s/n, 03202-Elche, Spain}\\{\small E.mail: jvalero@umh.es}}
\date{}
\maketitle

\begin{abstract}
This paper is devoted to the study of nonautonomous multivalued semiflows and
their associated pullback attractors. For this kind of dynamical systems we
are able to characterize the upper and lower bounds of the attractor as
complete trajectories belonging to the attractor, so that all the internal
dynamics is confined in this region, which can be described as an interval due
to the orderly nature of the processes. Thus, we are able to generalize to
this framework previous general results in literature for autonomous
multivalued flows or nonautonomous differential equations. We apply our
results to a partial differential inclusion with a nonautonomous term, also
proving the upper semicontinuity dependence of pullback and global attractors
when the time dependent term asymptotically converges to an autonomous
multivalued term.

\end{abstract}








\section{Introduction}

The study of attractors for infinite dimensional dynamical systems has been
one of the most active areas in dynamical system theory in the last five
decades. Nowadays the theory is very well suited with very good reference
texts (see, for instance, \cite{Ha, BV, He, La, Rob, SY, Vi92, Te, CV}). Among
the very interesting problems in this area, the one trying to characterize the
internal dynamics of attractors has been receiving a lot of interest (see
\cite{Henry85, Rob, CLRBook, CI, Ha, HR, HR92, HMO, ARV}), being the
topological and geometrical description of attractors the core subject in
order to understand the dynamics inside attracting compact sets (see
\cite{aragao, aragao2, BCCL, BCLBook, BCLR, CL09, CLRS, CLR12, conley, CV16}).
When the associated model of ordinary or partial differential equations is
autonomous, monotone systems is a class of associated semigroups for which a
first characterization of attractor can be described (see \cite{HS,
Chueshov}). Indeed, as we can define an order in the phase space, solutions
are also ordered and it can be proved that there exists an interval limited by
equilibria describing the attractor, i.e., there is no other point of the
attractor outside these upper and lower stationary solutions, and, moreover,
these solutions are asymptotically stable from above and below respectively.
These results have also been generalized to the nonautonomous and random
frameworks (see \cite{Chueshov, LangaSuarez2002, ArnoldChu,RRV,RBV}). Also, in
\cite{BrCarVal} the structure of the nonautonomous attractor for a scalar
one-dimensional parabolic equation is analysed in much more detail.

On the other hand, when uniqueness of solutions to a differential equation is
sometimes unknown, or equations include multivalued functionals, monotone
dynamical systems has also been developed for these cases, by generalizing the
single-valued results for multivalued operators or generalized semiflows (see
\cite{Valero2012, CLV2005}). We note that in this case, the way in which we
define an order for solutions is crucial, as different options are possible.
In this paper, we follow the approach given in \cite{CLV2005} to define
multivalued order preserving processes (see Definition \ref{def1}) in the
strongest sense. There are several previous papers describing the attractors
for nonautonomous multivalued flows in concrete examples of ordinary or
partial differential equations (see \cite{Valero2021, CLV2016, SV, SS,
CLV2020}). In these papers the concept of equilibria is generalized to the one
of complete trajectory, allowing for the definition of sub and super
trajectories (see Definition \ref{def2}). We take advantage of these concepts
for our main result (Theorem \ref{Structure}), for which we are able to give a
general result about existence of complete trajectories $\gamma_{\ast}\left(
\cdot\right)  ,\ \gamma^{\ast}(\cdot)$ such that $\gamma_{\ast}\left(
t\right)  ,\gamma^{\ast}(t)\in\mathcal{A}(t)$ and $\mathcal{A}(t)\subset
I_{\gamma_{\ast}}^{\gamma^{\ast}}\left(  t\right)  $ for all $t\in\mathbb{R}$,
being $I_{\gamma_{\ast}}^{\gamma^{\ast}}$ the interval generated by
$\gamma_{\ast}\left(  \cdot\right)  $ and $\gamma^{\ast}(\cdot)$.

We apply our general approach to study the problem
\begin{equation}
\left\{
\begin{array}
[c]{l}%
\dfrac{\partial u}{\partial t}-\dfrac{\partial^{2}u}{\partial x^{2}}\in
b(t)H_{0}(u)+\omega(t)u,\text{ on }(\tau,\infty)\times\left(  0,1\right)  ,\\
u(t,0)=u(t,1)=0,\\
u(\tau,x)=u_{\tau}(x),
\end{array}
\right.  \label{Incl0}%
\end{equation}
where $b:\mathbb{R}\rightarrow\mathbb{R}^{+},$ $\omega:\mathbb{R}%
\rightarrow\mathbb{R}^{+}$ are continuous functions satisfying suitable
conditions and $H_{0}(u)$ is the Heaviside function. We apply Theorem
\ref{Structure} to this problem (Section \ref{sec31}). Finally, in Section
\ref{sec32} we consider the situation when problem (\ref{Incl0}) is
asymptotically autonomous, proving the upper semicontinuity of pullback
attractors when times goes to $+\infty$.

\section{Nonautonomous monotone multivalued dynamical systems}

Let $X$ be a partially ordered complete metric space with metric $\rho$ and
let the order relation $\leq$ be compatible with the topology in the following
sense:\newline(1) For any bounded set B, there exist $a,b\in X$ such that
$B\subset\lbrack a,b]:=\{x\in X;\;a\leq x\leq b\};$\newline(2) If
$x_{n}\rightarrow x,\;y_{n}\rightarrow y$ and $x_{n}\leq y_{n},$ then $x\leq
y;$\newline(3) If $u\leq v\leq w,$ then $\rho(u,v)\leq\rho(u,w)$ and
$\rho(v,w)\leq\rho(u,w).$

Denote by $P(X)$ the set of all non-empty subsets of $X$ and let
$\mathbb{R}_{d}=\{(t,s)\in\mathbb{R}^{2}:t\geq s\}$. The Hausdorff
semidistance from the set $B$ to the set $A$ is given by $dist_{X}\left(
B,A\right)  =\sup_{b\in B}\inf_{a\in A}\rho\left(  b,a\right)  .$

We recall that the map $U:\mathbb{R}_{d}\times X\rightarrow P(X)$ is a
multivalued process if:

\begin{enumerate}
\item $U(t,t,$\textperiodcentered$)$ is the identity map for all
$t\in\mathbb{R};$

\item $U(t,s,x)\subset U(t,\tau,U(\tau,s,x))$ for any $x\in X,$ $s\leq\tau\leq
t.$
\end{enumerate}

Here, for any subset $B\subset X$, $U(t,\tau,B)$ stands for $\cup_{x\in
B}U(t,\tau,x)$. Also, $U$ is said to be a strict multivalued process if,
moreover,
\[
U(t,s,x)=U(t,\tau,U(\tau,s,x))\text{ for any }x\in X,\ s\leq\tau\leq t.
\]

In applications, $U$ is usually generated by the trajectories of a
differential equation. Such situation can be described by the following
abstract setting. Let $W_{\tau}=C([\tau,+\infty),X)$ and let $\mathcal{R}%
=\{\mathcal{R}_{\tau}\}_{\tau\in\mathbb{R}}$, $\mathcal{R}_{\tau}\subset
W_{\tau}$, be a family of functions. We consider the following axioms:

\begin{enumerate}
\item[$\left(  K1\right)  $] For any $\tau\in\mathbb{R}$ and $x\in X$ there
exists $\varphi\in\mathcal{R}_{\tau}$ such that $\varphi\left(  \tau\right)
=x.$

\item[$\left(  K2\right)  $] $\varphi^{+s}:=\varphi\mid_{\lbrack\tau
+s,\infty)}\in\mathcal{R}_{\tau+s}$ for any $s\geq0$, $\varphi\in
\mathcal{R}_{\tau}$ (translation property).

\item[$\left(  K3\right)  $] Let $\varphi,\psi\in\mathcal{R}$ be such that
$\varphi\in\mathcal{R}_{\tau}$, $\psi\in\mathcal{R}_{r}$ and $\varphi
(s)=\psi(s)$ for some $s\geq r\geq\tau$. Then the function $\theta$ defined
by
\[
\theta(t):=\left\{
\begin{array}
[c]{l}%
\varphi(t),\ t\in\lbrack\tau,s],\\
\psi(t),\ t\in\lbrack s,+\infty),
\end{array}
\right.
\]
belongs to $\mathcal{R}_{\tau}$ (concatenation property).

\item[$\left(  K4\right)  $] For any sequence $\varphi^{n}\in\mathcal{R}%
_{\tau}$ such that $\varphi^{n}\left(  \tau\right)  \rightarrow\varphi_{0} $
in $X$ as $n\to\infty$, there exists a subsequence $\varphi^{n_{k}}$ and
$\varphi\in\mathcal{R}_{\tau}$ such that $\varphi(\tau)=\varphi_{0}$ and
\[
\varphi^{n_{k}}\left(  t\right)  \rightarrow\varphi\left(  t\right)  \text{,
}\forall t\geq\tau.
\]

\end{enumerate}

We define the multivalued family of operators $U:\mathbb{R}_{d}\times
X\rightarrow P(X)$ associated with $\mathcal{R}$ by%
\[
U(t,s,x)=\{u(t):u(\text{\textperiodcentered})\in\mathcal{R}_{s},\text{
}u(s)=x\}.
\]
It follows from $\left(  K1\right)  -\left(  K2\right)  $ that $U$ is a
multivalued process. If, in addition, $(K3)$ is satisfied, then $U$ is strict.
Moreover, $\left(  K4\right)  $ implies that the graph of the map $x\mapsto
U(t,s,x)$ is closed for all $\left(  t,s\right)  \in\mathbb{R}_{d}$ .

We need to generalize to multivalued processes the concept of order preserving
process \cite{LangaSuarez2002}. For a single-valued process $S:\mathbb{R}%
_{d}\times X\rightarrow X$ this means that $S(t,s,x)\leq S(t,s,y)$ for all
$t\geq s$ and $y\geq x$. That is, the order of the initial datum is preserved
for all future times. However, in the multivalued case, when uniqueness of
solutions is not guaranteed, such definition is more complicated to establish,
as different concepts for comparison of solutions are possible (see
\cite{Valero2012}). In order to get results about the structure of pullback
attractor we need the strongest version of comparison, given in the following definition.

\begin{definition}
\label{def1} Assume that $\left(  K1\right)  -\left(  K2\right)  $ hold. The
multivalued process $U$ is said to be strongly order preserving if for any
initial data $x_{\tau}\leq y_{\tau}$ and $\tau\in\mathbb{R}$ there exist
\underline{$\varphi$}$,\overline{\varphi}\in\mathcal{R}_{\tau}$ such that
\underline{$\varphi$}$\left(  \tau\right)  =x_{\tau},\ \overline{\varphi
}\left(  \tau\right)  =y_{\tau}$ and%
\[
\underline{\varphi}(t)\leq y\left(  t\right)  \text{,\ }x\left(  t\right)
\leq\overline{\varphi}(t),\ \text{for all }t\geq\tau,
\]
where $x\left(  \text{\textperiodcentered}\right)  ,\ y\left(
\text{\textperiodcentered}\right)  \in\mathcal{R}_{\tau}$ satisfy $x\left(
\tau\right)  =x_{\tau}$, $y\left(  \tau\right)  =y_{\tau}$ and are arbitrary.
\end{definition}

This definition obviously implies that for any $x_{\tau}\leq y_{\tau}$ and
$\tau\leq t$ there exist elements $\overline{\varphi}\left(  t\right)  \in
U\left(  t,\tau,y_{\tau}\right)  ,\ \underline{\varphi}\left(  t\right)  \in
U(t,\tau,x_{\tau})$ such that%
\[
\underline{\varphi}(t)\leq y\text{,\ }x\leq\overline{\varphi}(t),\text{ }%
\]
for all $y\in U\left(  t,\tau,y_{\tau}\right)  ,\ x\in U(t,\tau,x_{\tau}).$
This could be the definition of order-preserving process for a general
multivalued process $U$, in a similar way as in the autononomous case
\cite{CLV2005}. However, when studying the structure of pullback attractors,
we work with complete trajectories, so this stronger definition is mandatory.

Prior to stating and proving our main result in this section, we will recall
some definitions and known results concerning pullback attractors for
multivalued processes.

The family of compact sets $\{\mathcal{A}(t)\}_{t\in\mathbb{R}}$ is called a
pullback attractor for $U$ if:

\begin{enumerate}
\item It is pullback attracting, which means that
\[
dist_{X}(U(t,s,B),\mathcal{A}(t))\rightarrow0,\text{ as }s\rightarrow-\infty,
\]
for any bounded set $B.$

\item $\mathcal{A}(t)\subset U(t,s,\mathcal{A}(s))$ for all $s\leq t$
(negative invariance).

\item $\{\mathcal{A}(t)\}_{t\in\mathbb{R}}$ is the minimal pullback attracting
family, that is, if $\{\mathcal{K}(t)\}_{t\in\mathbb{R}}$ is another pullback
attracting family of closed sets, then $\mathcal{A}(t)\subset\mathcal{K}(t)$
for any $t\in\mathbb{R}.$
\end{enumerate}

The pullback attractor is invariant if, moreover, $\mathcal{A}%
(t)=U(t,s,\mathcal{A}(s))$ for all $s\leq t.$ If the pullback attractor
$\{\mathcal{A}(t)\}_{t\in\mathbb{R}}$ is backwards bounded, that is,
$\cup_{s\leq\tau}\mathcal{A}(s)$ is a bounded set for some $\tau\in\mathbb{R}%
$, and the multivalued process $U$ is strict, then $\{\mathcal{A}%
(t)\}_{t\in\mathbb{R}} $ is invariant \cite[Lemma 2.5]{CLV2016}.

Further, we will characterize the pullback attractor in terms of bounded
complete trajectories.

The map $\gamma:\mathbb{R}\rightarrow X$ is called a complete orbit (sometimes
also called complete trajectory) of $\mathcal{R}$ if
\[
\varphi=\gamma\mid_{\lbrack\tau,+\infty)}\in\mathcal{R}_{\tau}\text{ for all
}\tau\in\mathbb{R}.
\]
A complete orbit $\gamma$ is bounded if $\cup_{t\in\mathbb{R}}\gamma\left(
t\right)  $ is a bounded set.

\begin{lemma}
\label{characterization} \cite[Corollaries 2.10 and 2.12]{CLV2016} Let
$\left(  K1\right)  -\left(  K2\right)  $ hold and that either $\left(
K3\right)  $ or $\left(  K4\right)  $ be satisfied. Assume that $U$ possesses
a pullback attractor $\{\mathcal{A}(t)\}_{t\in\mathbb{R}}$ which is bounded,
that is, $\cup_{t\in\mathbb{R}}\mathcal{A}(t)$ is a bounded set. Then
\[
\mathcal{A}(t)=\{\gamma\left(  t\right)  :\gamma\text{ is a bounded complete
trajectory}\}.
\]

\end{lemma}

General assumptions ensuring the existence of pullback attractors for
multivalued processes can be found for example in \cite{CLMV2003},
\cite{KapKasVal11}. See also \cite{SV} for more results on characterization of
pullback attractors. Assuming that a pullback attractor exists, we are
interested in proving the existence of two special complete orbits which are
upper and lower bounds of the pullback attractor, giving some insight into its
structure. In this way, we generalize to the multivalued case some well-known
results for single-valued processes \cite{Chueshov}, \cite{LangaSuarez2002}.

\begin{definition}
\label{def2} We say that the function \underline{$u$}$:\mathbb{R}\rightarrow
X$ is a complete sub-trajectory for $U$ if%
\[
\underline{u}\left(  t\right)  \leq y\text{ }\forall y\in U(t,s,\underline{u}%
\left(  s\right)  ),\ \forall t\geq s.
\]
The function $\overline{u}:\mathbb{R}\rightarrow X$ is a complete
super-trajectory for $U$ if%
\[
\overline{u}\left(  t\right)  \geq y\text{ }\forall y\in U(t,s,\overline
{u}\left(  s\right)  ),\ \forall t\geq s.
\]

\end{definition}

If \underline{$u$},\ $\overline{u}$ are sub and upper complete trajectories
for $U$ such that
\begin{equation}
\underline{u}\left(  t\right)  \leq\overline{u}\left(  t\right)  ,\text{ for
all }t\in\mathbb{R}, \label{OrderSubUpper}%
\end{equation}
$\ $we can define the interval%
\begin{equation}
I_{\underline{u}}^{\overline{u}}\left(  t\right)  =\{y\in X:\underline{u}%
\left(  t\right)  \leq y\leq\overline{u}\left(  t\right)  \}. \label{Interval}%
\end{equation}

\begin{lemma}
Assume that $\left(  K1\right)  -\left(  K2\right)  $ hold and that the
multivalued process $U$ is strongly order-preserving. Let \underline{$u$%
},\ $\overline{u}$ be sub and upper complete trajectories for $U$ satisfying
(\ref{OrderSubUpper}). Then the interval $I_{\underline{u}}^{\overline{u}}$ is
positively invariant, that is,%
\[
U(t,s,I_{\underline{u}}^{\overline{u}}(s))\subset I_{\underline{u}}%
^{\overline{u}}(t)\text{ for all }t\geq s.
\]

\end{lemma}

\begin{proof}
We fix any $s\in\mathbb{R}.$ Let $y\in U(t,s,I_{\underline{u}}^{\overline{u}%
}(s))$, $t\geq s$ be arbitrary. Then there exist $y_{s}\in I_{\underline{u}%
}^{\overline{u}}(s)$ and $\varphi\in\mathcal{R}_{s}$ such that $\varphi\left(
s\right)  =y_{s}$ and $y=\varphi\left(  t\right)  $. Since $U$ is strongly
order preserving and \underline{$u$},\ $\overline{u}$ are sub and upper
complete trajectories, there are \underline{$\varphi$}$,\overline{\varphi}%
\in\mathcal{R}_{s}$ such that \underline{$\varphi$}$\left(  s\right)
=\underline{u}\left(  s\right)  ,\ \overline{\varphi}\left(  s\right)
=\overline{u}\left(  s\right)  $ and%
\[
\underline{u}\left(  t\right)  \leq\underline{\varphi}(t)\leq y\leq
\overline{\varphi}(t)\leq\overline{u}\left(  t\right)  .
\]
Thus, $y\in I_{\underline{u}}^{\overline{u}}(t).$
\end{proof}

\bigskip

We are now ready to prove the main theorem of this work.

\begin{theorem}
\label{Structure}Assume that $\left(  K1\right)  ,\;\left(  K2\right)  $ and
$\left(  K4\right)  $ hold and that the multivalued process $U$ is strongly
order preserving and possesses a pullback attractor $\{\mathcal{A}%
(t)\}_{t\in\mathbb{R}}$. Let \underline{$u$},\ $\overline{u}$ be sub and upper
complete trajectories for $U$ satisfying (\ref{OrderSubUpper}) and such that
for all $t\in\mathbb{R},$%
\begin{equation}
\mathcal{A}(t)\subset I_{\underline{u}}^{\overline{u}}(t),
\label{AttrInterval}%
\end{equation}%
\begin{align}
dist_{X}\left(  U(t,s,\underline{u}\left(  s\right)  ),\mathcal{A}(t)\right)
&  \rightarrow0\text{,}\nonumber\\
dist_{X}\left(  U(t,s,\overline{u}\left(  s\right)  ),\mathcal{A}(t)\right)
&  \rightarrow0\text{, as }s\rightarrow-\infty. \label{AttrSubUpper}%
\end{align}
Then there exist complete orbits $\gamma_{\ast}\left(  \cdot\right)
,\gamma^{\ast}(\cdot)$ such that:

\begin{enumerate}
\item $\gamma_{\ast}\left(  t\right)  ,\gamma^{\ast}(t)\in\mathcal{A}(t)$ and
$\mathcal{A}(t)\subset I_{\gamma_{\ast}}^{\gamma^{\ast}}\left(  t\right)  $
for all $t\in\mathbb{R}.$ If, moreover, the pullback attractor is bounded,
i.e., if $\cup_{s\in\mathbb{R}}\mathcal{A}(s)$ is a bounded set in $X$, then
$\gamma_{\ast}(\cdot)$ and $\gamma^{\ast}(\cdot)$ are bounded.

\item \underline{$u$}$\left(  t\right)  \leq\gamma_{\ast}\left(  t\right)
\leq\gamma^{\ast}\left(  t\right)  \leq$\ $\overline{u}\left(  t\right)  $ for
any $t\in\mathbb{R}.$

\item $\gamma_{\ast}$ ($\gamma^{\ast}$) is minimal (maximal) in the sense that
there is no other complete orbit $\gamma$ such that
\[
\underline{u}\left(  t\right)  \leq\gamma\left(  t\right)  \leq\gamma_{\ast
}\left(  t\right)  \text{ (}\gamma^{\ast}\left(  t\right)  \leq\gamma\left(
t\right)  \leq\ \overline{u}\left(  t\right)  \text{) for all }t\in
\mathbb{R}\text{.}%
\]

\item Assume that $\varphi_{s}:=\gamma_{\ast}\mid_{\lbrack s,+\infty)}$
($\gamma^{\ast}\mid_{\lbrack s,+\infty)}$) is the unique function in
$\mathcal{R}_{s}$ such that $\varphi_{s}\left(  s\right)  =\gamma_{\ast}(s)$
($\gamma^{\ast}(s)$) for any $s\in\mathbb{R}$. Then $\gamma_{\ast}$ is
globally asymptotically stable from below, i.e.,
\begin{equation}
\lim_{\tau\rightarrow-\infty}dist(U(t,\tau,v_{\tau}),\gamma_{\ast}(t))=0
\label{StableBelow}%
\end{equation}
whenever $v_{\tau}\in X$ with $\underline{u}(\tau)\leq v_{\tau}\leq
\gamma_{\ast}(\tau)\;\forall\;\tau\in\mathbb{R};$ and $\gamma^{\ast}$ is
globally asymptotically stable from above, i.e.,
\begin{equation}
\lim_{\tau\rightarrow-\infty}dist(U(t,\tau,v_{\tau}),\gamma^{\ast}(t))=0
\label{StableAbove}%
\end{equation}
whenever $v_{\tau}\in X$ with $\overline{u}(\tau)\geq v_{\tau}\geq\gamma
^{\ast}(\tau)\;\forall\;\tau\in\mathbb{R}.$
\end{enumerate}
\end{theorem}

\begin{proof}
\textit{Proof of item 1:} Let us consider a sequence $\{s_{n}\}_{n\in
\mathbb{N}}$ such that $s_{n}\rightarrow-\infty$ as $n\rightarrow+\infty. $
From (\ref{AttrInterval}), $\underline{u}(s_{n})\leq u\leq\overline{u}%
(s_{n}),$ for all $u\in\mathcal{A}(s_{n}),\;n\in\mathbb{N}.$ Since $U$ is
strongly order preserving, there exist solutions $\varphi_{n}\in
\mathcal{R}_{s_{n}},$ $\varphi^{n}\in\mathcal{R}_{s_{n}}$ such that
$\varphi_{n}(s_{n})=\underline{u}(s_{n}),\ \varphi^{n}\left(  s_{n}\right)
=\overline{u}(s_{n})$ and%
\begin{equation}
\varphi_{n}(t)\leq y\left(  t\right)  \leq\varphi^{n}(t)\text{ }
\label{IneqFin}%
\end{equation}
for any $t\geq s_{n}$ and $y\in\mathcal{R}_{s_{n}}$ satisfying $\underline{u}%
(s_{n})\leq y\left(  s_{n}\right)  \leq\overline{u}(s_{n})$. The negative
invariance of the pullback attractor implies then that%
\begin{equation}
\varphi_{n}(t)\leq\xi\leq\varphi^{n}(t),\;\forall\;\xi\in\mathcal{A}%
(t),\;t\geq s_{n}\text{, }n\in\mathbb{N}. \label{qaz}%
\end{equation}

Let us obtain the orbit $\gamma_{\ast}\left(  \text{\textperiodcentered
}\right)  $. Since $\mathcal{A}(t)$ is a compact set, by using
(\ref{AttrSubUpper}), we can choose converging subsequences, which are denoted
the same, such that $\varphi_{n}(0)\rightarrow z_{\ast}^{0}\in\mathcal{A}(0)$.
As $\varphi_{n}\in\mathcal{R}_{s_{n}}$, from the property $(K2)$,
\[
\varphi_{n}^{+(-s_{n})}=\varphi_{n}\mid_{\lbrack0,\infty)}\in\mathcal{R}%
_{0},\;\forall\;n\in\mathbb{N}.
\]
By the property $(K4)$ there exist a subsequence of $\{\varphi_{n}^{+(-s_{n}%
)}\}_{n\in\mathbb{N}},$ which we do not relabel, and a solution $\psi_{0}%
\in\mathcal{R}_{0}$ with $\psi_{0}(0)=z_{\ast}^{0}$ such that
\[
\varphi_{n}^{+(-s_{n})}(t)=\varphi_{n}(t)\rightarrow\psi_{0}(t)\in
\mathcal{A}(t),\;\forall\;t\geq0.
\]
Passing the limit in (\ref{qaz}) and using that the order relation is
compatible with the topology of the space, we obtain
\[
\psi_{0}\left(  t\right)  \leq u,\;\forall\;u\in\mathcal{A}(t),\ t\geq0.
\]
Further, passing to a new subsequence we can state that $\varphi
_{n}(-1)\rightarrow z_{\ast}^{-1}\in\mathcal{A}(-1)$. Again, by $\left(
K2\right)  $ and $\left(  K4\right)  ,$ we have%
\[
\varphi_{n}^{+(-1-s_{n})}=\varphi_{n}\mid_{\lbrack-1,\infty)}\in
\mathcal{R}_{-1},\;\forall\;n\in\mathbb{N},
\]%
\[
\varphi_{n}^{+(-1-s_{n})}(t)=\varphi_{n}(t)\rightarrow\psi_{-1}(t)\in
\mathcal{A}(t),\;\forall\;t\geq-1,
\]
where $\psi_{-1}\in\mathcal{R}_{-1}$ with $\psi_{-1}(-1)=z_{\ast}^{-1}$. Also,
it is clear that%
\[
\psi_{-1}\left(  t\right)  \leq u,\;\forall\;u\in\mathcal{A}(t),\ t\geq-1.
\]%
\[
\psi_{-1}\left(  t\right)  =\psi_{0}\left(  t\right)  ,\ \forall\ t\geq0.
\]
Continuing in this way we obtain a sequence $\psi_{-k}\in\mathcal{R}_{-k}$,
$k=0,1,...$, such that%
\[
\psi_{-k}(t)\in\mathcal{A}(t),\;\forall\;t\geq-k,
\]%
\[
\psi_{-k}\left(  t\right)  \leq u,\;\forall\;u\in\mathcal{A}(t),\ t\geq-k,
\]%
\[
\psi_{-k}\left(  t\right)  =\psi_{-k+1}\left(  t\right)  ,\ \forall
\ t\geq-k+1.
\]
We define $\gamma_{\ast}\left(  \text{\textperiodcentered}\right)  $ by taking
the common value of the functions $\psi_{-k}$ for any $t\in\mathbb{R}$. It is
obvious that $\gamma_{\ast}\left(  \text{\textperiodcentered}\right)  $ is a
complete orbit satisfying%
\[
\gamma_{\ast}\left(  t\right)  \in\mathcal{A}(t),\;\forall\;t\in\mathbb{R},
\]%
\begin{equation}
\gamma_{\ast}\left(  t\right)  \leq u,\;\forall\;u\in\mathcal{A}%
(t),\ t\in\mathbb{R}. \label{wsx}%
\end{equation}

In a similar way we obtain a complete orbit $\gamma^{\ast}\left(
\text{\textperiodcentered}\right)  $ satisfying
\[
\gamma^{\ast}\left(  t\right)  \in\mathcal{A}(t),\;\forall\;t\in\mathbb{R},
\]%
\begin{equation}
\gamma^{\ast}\left(  t\right)  \geq u,\;\forall\;u\in\mathcal{A}%
(t),\ t\in\mathbb{R}. \label{wsx2}%
\end{equation}

Finally, we observe that by a diagonal argument there is a subsequence such
that%
\begin{align}
\varphi_{n}(t)  &  \rightarrow\gamma_{\ast}\left(  t\right)
,\label{ConvergFin}\\
\varphi^{n}(t)  &  \rightarrow\gamma^{\ast}\left(  t\right)  ,\text{ }\forall
t\in\mathbb{R}.\nonumber
\end{align}

\textit{Proof of item 2:} We already know from the previous item (see
(\ref{wsx})) that $\gamma_{\ast}(t)\leq\gamma^{\ast}(t).$ Let us prove that
$\underline{u}\left(  t\right)  \leq\gamma_{\ast}\left(  t\right)  $ (the
proof of $\gamma^{\ast}\left(  t\right)  \leq\overline{u}\left(  t\right)  $
is analogous). Since $\underline{u}$ is a complete sub-trajectory for $U$, we
have $\underline{u}\left(  t\right)  \leq y\text{ for all }y\in
U(t,s,\underline{u}\left(  s\right)  ),\ t\geq s.$ Fix $s<t$. By using
hypothesis (\ref{AttrInterval}), we have $\underline{u}(s)\leq v\leq
\overline{u}(s),$ for all $v\in\mathcal{A}(s).$ As $U$ is strongly order
preserving, there exists $y_{\ast}\in U(t,s,\underline{u}\left(  s\right)  )$
such that $y_{\ast}\leq\gamma(t)$ for all $\gamma(t)\in U(t,s,v)$ and any
$v\in\mathcal{A}(s).$ In particular, as $\gamma_{\ast}\mid_{\lbrack s,\infty
)}\in\mathcal{R}_{s}$ and $\gamma_{\ast}\left(  s\right)  \in\mathcal{A}(s)$,
we have
\[
\underline{u}\left(  t\right)  \leq y_{\ast}\leq\gamma_{\ast}(t).
\]

\textit{Proof of item 3:} Let $\gamma(\cdot)$ be a complete orbit such that
\[
\underline{u}\left(  t\right)  \leq\gamma\left(  t\right)  \leq\gamma_{\ast
}\left(  t\right)  \;\text{for all}\;t\in\mathbb{R}.
\]
By (\ref{IneqFin}), we infer that%
\[
\varphi_{n}(t)\leq\gamma(t)\ \ \forall t\geq s_{n}.
\]
Then, by using (\ref{ConvergFin}), we get%
\[
\gamma_{\ast}\left(  t\right)  =\lim_{n\rightarrow+\infty}\varphi_{n}%
(t)\leq\gamma\left(  t\right)  \leq\gamma_{\ast}\left(  t\right)  .
\]
Therefore, $\gamma_{\ast}\left(  t\right)  =\gamma\left(  t\right)  $ for all
$t\in\mathbb{R}.$ Analogously we can prove that there is no other complete
orbit $\gamma$ such that $\gamma^{\ast}\left(  t\right)  \leq\gamma\left(
t\right)  \leq\ \overline{u}\left(  t\right)  \text{ for all }t\in\mathbb{R}.$

\textit{Proof of item 4:} Consider $\{v_{s_{n}}\}_{n\in\mathbb{N}}\subset X$
with $\underline{u}(s_{n})\leq v_{s_{n}}\leq\gamma_{\ast}(s_{n})$ for each
$n\in\mathbb{N}$ and $s_{n}\rightarrow-\infty$ as $n\rightarrow\infty.$ By
(\ref{IneqFin}) we have%
\[
\varphi_{n}(t)\leq y\text{ for all }y\in U(t,s_{n},v_{s_{n}}).
\]
Also, since $U$ is strongly order preserving and $\varphi_{s_{n}}=\gamma
_{\ast}\mid_{\lbrack s_{n},+\infty)}$ is the unique function in $\mathcal{R}%
_{s_{n}}$ with initial value $\gamma_{\ast}(s_{n})$, we obtain that%
\begin{equation}
y\leq\gamma_{\ast}(t)\text{ for all }y\in U(t,s_{n},v_{s_{n}}).
\label{gammaAbove}%
\end{equation}
Then by the third condition of compatibility of the order with the topology,
we have that $\rho(y,\gamma_{\ast}(t))\leq\rho(\varphi_{n}(t),\gamma_{\ast
}(t))\;$for all $n\in\mathbb{N}.$ Thus,
\[
dist(U(t,s_{n},v_{s_{n}}),\gamma_{\ast}(t))=\sup_{y\in U(t,s_{n},v_{\tau_{n}%
})}\rho(y,\gamma_{\ast}(t))\leq\rho(\varphi_{n}(t),\gamma_{\ast}%
(t))\rightarrow0
\]
as $n\rightarrow\infty.$ Therefore, $\gamma_{\ast}$ is globally asymptotically
stable from below. The proof of the fact that $\gamma^{\ast}$ be globally
asymptotically stable from above is analogous and we leave it to the reader.
\end{proof}

\bigskip

\begin{remark}
We can avoid using the third condition of compatibility in the proof of item 4
if we relax a bit the definition of $\gamma_{\ast}$ ($\gamma^{\ast}$) being
globally asymptotically stable from below (above). Namely, if
(\ref{StableBelow}) has to be satisfied only whenever $v_{\tau}\in B$, where
$B$ is a bounded set, then assuming that (\ref{StableBelow}) is not true there
would exist $\varepsilon>0$ and $y_{n}\in U(t,s_{n},v_{s_{n}})$, where
$s_{n}\rightarrow-\infty$ as $n\rightarrow\infty$, such that%
\begin{equation}
\rho(y_{n},\gamma_{\ast}(t))\geq\varepsilon\text{ for all }n\geq n_{0}\text{.}
\label{Contr}%
\end{equation}
By (\ref{IneqFin}) and (\ref{gammaAbove}) we have%
\[
\varphi_{n}(t)\leq y\leq\gamma_{\ast}(t)\text{ for all }y\in U(t,s_{n}%
,v_{s_{n}}).
\]
By the attracting property of the pullback attractor, passing to a subsequence
we have that $y_{n}\rightarrow y_{0}\in\mathcal{A}(t).$ Hence, the second
condition of compatibility gives%
\[
y_{0}\leq\gamma_{\ast}(t)\leq y_{0},
\]
so $\gamma_{\ast}(t)=y_{0}$, which contradicts (\ref{Contr}). The same
argument is valid for $\gamma^{\ast}.$
\end{remark}

\bigskip

We finish this section by generalizing a result from \cite{SS} (see Theorem
4.1) about upper semicontinuity of pullback attractors for asymptotically
autonomous systems.

We recall that for a family of functions $\mathcal{R}_{0}\subset
C([0,+\infty),X)$ axioms $\left(  K1\right)  -\left(  K4\right)  $ read in the
autonomous situation as follows:

\begin{enumerate}
\item[$\left(  H1\right)  $] For any $x\in X$ there exists $\varphi
\in\mathcal{R}_{0}$ such that $\varphi\left(  0\right)  =x.$

\item[$\left(  H2\right)  $] $\varphi^{+s}:=\varphi\left(
\text{\textperiodcentered}+s\right)  \in\mathcal{R}_{0}$ for any $s\geq0$,
$\varphi\in\mathcal{R}_{0}$ (translation property).

\item[$\left(  H3\right)  $] Let $\varphi,\psi\in\mathcal{R}_{0}$ be such that
$\varphi(s)=\psi(0)$ for some $s>0$. Then the function $\theta$ defined by
\[
\theta(t):=\left\{
\begin{array}
[c]{l}%
\varphi(t),\ t\in\lbrack0,s],\\
\psi(t-s),\ t\in\lbrack s,+\infty),
\end{array}
\right.
\]
belongs to $\mathcal{R}_{0}$ (concatenation property).

\item[$\left(  H4\right)  $] For any sequence $\varphi^{n}\in\mathcal{R}_{0}$
such that $\varphi^{n}\left(  0\right)  \rightarrow\varphi_{0}$ in $X$ as
$n\rightarrow\infty$, there exists a subsequence $\varphi^{n_{k}}$ and
$\varphi\in\mathcal{R}_{0}$ such that $\varphi(0)=\varphi_{0}$ and
\[
\varphi^{n_{k}}\left(  t\right)  \rightarrow\varphi\left(  t\right)  \text{,
}\forall t\geq0.
\]

\end{enumerate}

We define the multivalued family of operators $G:\mathbb{R}^{+}\times
X\rightarrow P(X)$ associated with $\mathcal{R}_{0}$ by%
\[
G(t,x)=\{u(t):u(\text{\textperiodcentered})\in\mathcal{R}_{0},\text{
}u(0)=x\}.
\]
It follows from $\left(  H1\right)  -\left(  H2\right)  $ that $G$ is a
multivalued semiflow, that is, $G(0,x)=x,\ G(t+s,x)\subset G(t,G(s,x))$ for
all $x\in X$, $t,s\geq0$. If, in addition, $(H3)$ is satisfied, then $G$ is
strict, which means that, moreover, $G(t+s,x)=G(t,G(s,x))$ for all $x\in X$,
$t,s\geq0$. In addition, $\left(  H4\right)  $ implies that the graph of the
map $x\mapsto G(t,x)$ is closed for all $t\geq0$ .

A compact set $\mathcal{A}$ is said to be a global attractor for $G$ if:

\begin{enumerate}
\item (Attraction property) For any bounded set $B$ we have:%
\[
dist(G(t,B),\mathcal{A})\rightarrow0\text{ as }t\rightarrow+\infty.
\]

\item (Negative invariance) $\mathcal{A}\subset G(t,\mathcal{A})$ for all
$t\geq0.$
\end{enumerate}

The attractor is called invariant if $\mathcal{A}=G(t,\mathcal{A})$ for all
$t\geq0.$

The map $\gamma:\mathbb{R}\rightarrow X$ is called a complete orbit of
$\mathcal{R}_{0}$ if
\[
\varphi\left(  \text{\textperiodcentered}\right)  =\gamma\left(
\text{\textperiodcentered}+\tau\right)  \in\mathcal{R}_{0}\text{ for all }%
\tau\in\mathbb{R}.
\]
A complete orbit $\gamma$ is bounded if $\cup_{t\in\mathbb{R}}\gamma\left(
t\right)  $ is a bounded set. If $\left(  H1\right)  -\left(  H2\right)  $ and
either $\left(  H3\right)  $ or $\left(  H4\right)  $ hold, then it is known
\cite[Theorems 9 and 10]{KapKasVal14} that
\[
\mathcal{A}=\{\gamma\left(  0\right)  :\gamma\text{ is a bounded complete
orbit of }\mathcal{R}_{0}\}
\]
whenever the global attractor exists.

\begin{theorem}
\label{abstract_upper_semicont XX} Let $\mathcal{R}=\{\mathcal{R}_{\tau
}\}_{\tau\in\mathbb{R}}$ be a family satisfying $\left(  H1\right)  -\left(
H2\right)  $ and such that the associated semiprocess $U$ has a pullback
attractor $\{\mathcal{A}(t)\}_{t\in\mathbb{R}}$. Assume that $\mathcal{R}$ is
asymptotically autonomous in the sense that there exists a family
$\mathcal{R}_{0}$ satisfying $\left(  H1\right)  -\left(  H2\right)  $ such
that for every sequences $\tau_{n}\rightarrow+\infty,\ \varphi^{n}%
\in\mathcal{R}_{\tau_{n}}$ such that $\varphi^{n}\left(  \tau_{n}\right)
\in\mathcal{A}(\tau_{n})$ and $\varphi^{n}\left(  \tau_{n}\right)
\rightarrow\varphi_{0}$, there exists a subsequence of $\{\psi_{n}\left(
\text{\textperiodcentered}\right)  \},\ \psi_{n}\left(
\text{\textperiodcentered}\right)  :=\varphi^{n}\left(
\text{\textperiodcentered}+\tau_{n}\right)  $, and $\psi_{0}\in\mathcal{R}%
_{0}$ for which $\psi_{n}(t)\rightarrow\psi_{0}(t)$ for all $t\geq0$. If
$\overline{\cup_{s\geq\tau}\mathcal{A}(s)}$ is a compact set for all $\tau
\in\mathbb{R}$ and the semiflow $G$ generated by $\mathcal{R}_{0}$ has a
global attractor $\mathcal{A}$, then%
\begin{equation}
\lim_{t\rightarrow+\infty}dist(\mathcal{A}(t),\mathcal{A})=0.
\label{ConvergAttr}%
\end{equation}

\end{theorem}

\begin{proof}
We know by assumption that $K=\overline{\cup_{s\geq0}\mathcal{A}(s)}$ is
compact. If (\ref{ConvergAttr}) were not true, there would exist
$\varepsilon>0$ and sequences $s_{n}\nearrow+\infty,$ $a_{n}\in\mathcal{A}%
(s_{n})$ such that%
\begin{equation}
dist(a_{n},\mathcal{A})\geq3\varepsilon\text{\ for all }n.
\label{Contradiction}%
\end{equation}
Since $\mathcal{A}$ attracts the bounded set $K$, there exists $n_{0}$ such
that
\[
dist(G(s_{n},K),\mathcal{A})\leq\varepsilon\text{ for all }n\geq n_{0}.
\]
By the negative invariance of the pullback attractor, for each $n\geq n_{0}$
there exist $b_{n}\in\mathcal{A}(s_{n}-s_{n_{0}})\subset K$ such that%
\[
a_{n}\in U(s_{n},s_{n}-s_{n_{0}},b_{n}).
\]
By the compactness of $K$ up to a subsequence we have that $b_{n}%
\rightarrow\varphi_{0}\in K$. Let $\varphi^{n}\in\mathcal{R}_{\tau_{n}}$,
$\tau_{n}=s_{n}-s_{n_{0}}$, be such that $\varphi^{n}(\tau_{n})=b_{n}$ and
$\varphi^{n}(s_{n})=a_{n}$. Since $\mathcal{R}$ is asymptotically autonomous,
there exists $\psi_{0}\in\mathcal{R}_{0}$ such that $\psi_{0}(0)=\varphi_{0}$
and up to a subsequence $\varphi^{n}\left(  t+\tau_{n}\right)  \rightarrow
\psi_{0}\left(  t\right)  $ for all $t\geq0$. In particular, we obtain that%
\[
a_{n}=\varphi^{n}(s_{n})=\varphi^{n}\left(  s_{n_{0}}+\tau_{n}\right)
\rightarrow\psi_{0}\left(  s_{n_{0}}\right)  .
\]
Hence, there is $n_{1}\geq n_{0}$ such that $\rho\left(  a_{n},\psi_{0}\left(
s_{n_{0}}\right)  \right)  \leq\varepsilon$ for all $n\geq n_{1}.$ Thus, as
$\psi_{0}\left(  s_{n_{0}}\right)  \in G(s_{n_{0}},\varphi_{0})\subset
G(s_{n_{0}},K)$, we have\
\[
dist(a_{n},\mathcal{A})\leq\rho(a_{n},\psi_{0}\left(  s_{n_{0}}\right)
)+dist(G(s_{n_{0}},K),\mathcal{A})\leq2\varepsilon\text{ for all }n\geq
n_{1},
\]
which contradicts (\ref{Contradiction}).
\end{proof}

\section{Application}

\subsection{An application of Theorem~\ref{Structure}}

\label{sec31}

We consider the problem
\begin{equation}
\left\{
\begin{array}
[c]{l}%
\dfrac{\partial u}{\partial t}-\dfrac{\partial^{2}u}{\partial x^{2}}\in
b(t)H_{0}(u)+\omega(t)u,\text{ on } (\tau,\infty)\times\left(  0,1\right)  ,\\
u(t,0)=u(t,1)=0,\\
u(\tau,x)=u_{\tau}(x),
\end{array}
\right.  \label{Incl}%
\end{equation}
where $b:\mathbb{R}\rightarrow\mathbb{R}^{+},$ $\omega:\mathbb{R}%
\rightarrow\mathbb{R}^{+}$ are continuous functions such that
\[
0<b_{0}\leq b\left(  t\right)  \leq b_{1},0\leq\omega_{0}\leq\omega\left(
t\right)  \leq\omega_{1}<\pi^{2},
\]
and
\[
H_{0}(u)=\left\{
\begin{array}
[c]{ll}%
-1, & \text{if }u<0,\\
\left[  -1,1\right]  , & \text{if }u=0,\\
1, & \text{if }u>0,
\end{array}
\right.
\]
is the Heaviside function. We note that $\pi^{2}$ is the first eigenvalue of
the operator $-\dfrac{\partial^{2}}{\partial x^{2}}$ in $H_{0}^{1}\left(
0,1\right)  .$

We consider the phase space $X=L^{2}\left(  0,1\right)  $ and its usual order
$\leq$ given by%
\[
u\leq v\Leftrightarrow u\left(  x\right)  \leq v\left(  x\right)  \text{ for
a.a. }x\in\left(  0,1\right)  .
\]

Let $A:D(A)\rightarrow H,\ D(A)=H^{2}(0,1)\cap H_{0}^{1}(\Omega),$ be the
operator $A=-\dfrac{d^{2}}{dx^{2}}$ with Dirichlet boundary conditions. We say
that the function $u\in C([\tau,+\infty),L^{2}\left(  0,1\right)  )$ is a
strong solution of (\ref{Incl}) if:

\begin{enumerate}
\item $u(\tau)=u_{\tau}$;

\item $u\left(  \text{\textperiodcentered}\right)  $ is absolutely continuous
on $[T_{1},T_{2}]$ for any $\tau<T_{1}<T_{2}$ and $u\left(  t\right)  \in D(A)
$ for a.a. $t\in\left(  T_{1},T_{2}\right)  ;$

\item There exists a function $f\in L_{loc}^{2}(\tau,+\infty;L^{2}\left(
0,1\right)  )$ such that $f(t,x)\in H_{0}(u(t,x))$, for a.e. $(t,x)\in
(\tau,+\infty)\times\left(  0,1\right)  $, and
\begin{equation}
\frac{du}{dt}(t)-Au(t)=b(t)f(t)+\omega(t)u(t),\text{ for a.e. }t\in
(\tau,+\infty), \label{EqSol}%
\end{equation}
where the equality is understood in the sense of space $L^{2}(0,1).$
\end{enumerate}

The set $\mathcal{R}_{\tau}$ is given by all the strong solutions to problem
(\ref{Incl}) with initial condition in $L^{2}\left(  0,1\right)  $ at initial
time $\tau$. Properties $\left(  K1\right)  -\left(  K4\right)  $ are
satisfied \cite{CLV2020}. It is known \cite[Theorem 5 and Lemma 6]{CLV2020}
that the associated multivalued process $U$ possesses an invariant pullback
attractor $\{\mathcal{A}(t)\}_{t\in\mathbb{R}}$ such that $\cup_{t\in
\mathbb{R}}\mathcal{A}(t)$ is bounded in $H_{0}^{1}\left(  0,1\right)  $,
$\overline{\cup_{t\in\mathbb{R}}\mathcal{A}(t)}$ is compact in $L^{2}\left(
0,1\right)  $ and
\[
\mathcal{A}(t)=\{\gamma\left(  t\right)  :\gamma\text{ is a bounded complete
trajectory}\}.
\]
Moreover, $U$ is strongly order preserving \cite[Theorem 2]{CLV2020}.

Let us consider the particular autonomous problem%
\begin{equation}
\left\{
\begin{array}
[c]{l}%
\dfrac{\partial u}{\partial t}-\dfrac{\partial^{2}u}{\partial x^{2}}\in
b_{1}H_{0}(u)+\omega_{1}u,\text{ on }\left(  0,1\right)  \times(\tau
,\infty),\\
u(t,0)=u(t,1)=0,\\
u(\tau,x)=u_{\tau}(x),
\end{array}
\right.  \label{Aut}%
\end{equation}
that is, when $b\left(  t\right)  \equiv b_{1},\ \omega\left(  t\right)
\equiv\omega_{1}$. The following comparison principle between solutions of
problems (\ref{Aut}) and (\ref{Incl}) was proved in \cite[Theorem 3]{CLV2020}.
We recall that a solution $u\left(  \text{\textperiodcentered}\right)  $ is
non-negative if $u\left(  t\right)  \geq0$ for all $t\geq\tau.$

\begin{lemma}
\label{Comparison}For any initial datum $u_{\tau}\geq0$ there exists a
non-negative strong solution$\ \overline{u}_{b_{1},\omega_{1}}\left(
\text{\textperiodcentered}\right)  $ to problem (\ref{Aut}) such that%
\begin{equation}
v\left(  t\right)  \leq\overline{u}_{b_{1},\omega_{1}}\left(  t\right)
\text{, }\forall t\geq\tau, \label{IneqNonautAut1}%
\end{equation}
where $v\left(  \text{\textperiodcentered}\right)  $ is an arbitrary strong
non-negative solution to (\ref{Incl})\ with $v\left(  \tau\right)  =u_{\tau}.$
\end{lemma}

\begin{corollary}
\label{Comparison2}If $v_{\tau}\leq u_{\tau}$, $u_{\tau}\geq0$, then
\begin{equation}
v\left(  t\right)  \leq\overline{u}_{b_{1},\omega_{1}}\left(  t\right)
\text{, }\forall t\geq\tau, \label{IneqNonautAut2}%
\end{equation}
where $v\left(  \text{\textperiodcentered}\right)  $ is an arbitrary strong
solution to (\ref{Incl})\ with $v\left(  \tau\right)  =v_{\tau}.$
\end{corollary}

\begin{proof}
Since $U$ is strongly order preserving, there exists a maximal solution
$\overline{v}\left(  \text{\textperiodcentered}\right)  $ to (\ref{Incl}%
)\ with $\overline{v}\left(  \tau\right)  =u_{\tau}$ such that $v\left(
t\right)  \leq\overline{v}\left(  t\right)  $ for all $t\geq\tau$ and any
strong solution $v\left(  \text{\textperiodcentered}\right)  $ to
(\ref{Incl})\ with $\overline{v}\left(  \tau\right)  =v_{\tau}$. Since
$u_{\tau}\geq0$, problem (\ref{Incl}) possesses at least one non-negative
strong solution \cite[Corollary 5]{CLV2020}, so $\overline{v}\left(
\text{\textperiodcentered}\right)  $ has to be non-negative. Hence, Lemma
\ref{Comparison} implies (\ref{IneqNonautAut2}).
\end{proof}

\bigskip

We recall other results from \cite{CLV2020}. In the autonomous case, the
multivalued process $U$ reduces to the multivalued semiflow $G:\mathbb{R}%
^{+}\times X\rightarrow P(X)$ given by%
\[
G(t,u_{0})=U(t,0,u_{0})=\{u(t):u(\text{\textperiodcentered})\in\mathcal{R}%
_{0},\text{ }u(0)=u_{0}\}.
\]
This map is a strict multivalued semiflow, that is, $G(0,$\textperiodcentered
$)$ is the identity map and $G(t+s,u_{0})=G(t,G(s,u_{0}))$ for any $t,s\geq0$,
$u_{0}\in X$. It possesses a global invariant attractor $\mathcal{A}$. There
exist for the problem (\ref{Aut}) an infinite (but countable) number of
stationary points, which are described in detail in \cite{ARV}. We are
interested now in two special stationary points. Namely, problem (\ref{Aut})
has one positive fixed point $v_{1}^{+}$ ($v_{1}^{+}\left(  x\right)  >0$ for
$x\in\left(  0,1\right)  $) and one negative fixed point $v_{1}^{-}=-v_{1}%
^{+}$, and the interval generated by them contains the global attractor:%
\[
v_{1}^{-}\leq y\leq v_{1}^{+}\text{ for all }y\in\mathcal{A}.
\]

We shall establish that $v_{1}^{+}$ ($v_{1}^{-})$ is an upper (sub) complete
trajectory for $U$ and that the interval $I_{v_{1}^{-}}^{v_{1}^{+}}$ satisfies
(\ref{AttrInterval}). We observe that the solution starting at $v_{1}^{+}$
($v_{1}^{-}$) is unique, which follows from its stability \cite[Theorem
6.3]{ARV}. Although this result is proved in \cite{ARV} for the particular
case $\omega_{1}=0$, it is true for any $0\leq\omega_{1}<\pi^{2}.$

\begin{lemma}
\label{upper}$v_{1}^{+}$ is a complete upper-trajectory for $U$. $v_{1}^{-}$
is a complete sub-trajectory for $U$.
\end{lemma}

\begin{proof}
Let $s$ be arbitrary. Since $u\left(  r\right)  =v_{1}^{+}$, $r\geq s$, is the
unique solution to problem (\ref{Aut}), by Corollary \ref{Comparison2} any
strong solution $v\left(  \text{\textperiodcentered}\right)  $ to problem
(\ref{Incl}) with $v\left(  s\right)  =v_{1}^{+}$ satisfies $v\left(
r\right)  \leq v_{1}^{+}$ for all $r\geq0.$ Hence, $v_{1}^{+}\geq y$ for any
$y\in U(t,s,v_{1}^{+}),\ t\geq s.$

For the second statement, as\ the map $H_{0}$ is odd, we infer that $w\left(
t\right)  =-v\left(  t\right)  $ is a strong solution to (\ref{Incl}) for any
strong solution $v\left(  \text{\textperiodcentered}\right)  $ with $v\left(
s\right)  =v_{1}^{-}=-v_{1}^{+}$. Therefore, since $w\left(  s\right)
=v_{1}^{+}$, we have by the previous result that $w\left(  r\right)  \leq
v_{1}^{+}$ for all $r\geq s$. Thus, $v\left(  r\right)  \geq-v_{1}^{+}$ for
all $r\geq s$, so $y\geq v_{1}^{-}$ for any $y\in U(t,s,v_{1}^{-}),\ t\geq s.
$
\end{proof}

\begin{theorem}
\label{Intervalo}$\mathcal{A}(t)\subset I_{v_{1}^{-}}^{v_{1}^{+}}%
=\{y:v_{1}^{-}\leq y\leq v_{1}^{+}\}$ for all $t\in\mathbb{R}.$
\end{theorem}

\begin{proof}
Suppose, by contradiction, that there are $t_{0}\in\mathbb{R}$ and
$y\in\mathcal{A}(t_{0})$ such that $y\not \leq v_{1}^{+}$. Since $y,v_{1}%
^{+}\in H_{0}^{1}\left(  0,1\right)  \subset C\left(  [0,1]\right)  $, there
has to exist an interval $[x_{0},x_{1}]\subset\left(  0,1\right)  $ such that
$y\left(  x\right)  >v_{1}^{+}\left(  x\right)  $ for any $x\in\lbrack
x_{0},x_{1}]$. Also, $y=\gamma\left(  t_{0}\right)  $, where $\gamma$ is a
bounded complete orbit. Since $\cup_{t\in\mathbb{R}}\mathcal{A}(t)$ is bounded
in $H_{0}^{1}\left(  0,1\right)  \subset C([0,1])$, there exists $c_{0}>0$
such that $\gamma\left(  s\right)  \leq u_{0}\equiv c_{0}$ for any
$s\in\mathbb{R}$. Hence, Corollary \ref{Comparison2} implies for any $s<t_{0}$
the existence of a non-negative strong solution $u_{s}\left(
\text{\textperiodcentered}\right)  $ to problem (\ref{Aut}) such that
$u_{s}\left(  s\right)  =u_{0}$ and $\gamma\left(  r\right)  \leq u_{s}\left(
r\right)  $ for all $r\geq s$.

The global attractor $\mathcal{A}$ attracts any bounded set of $L^{2}\left(
0,1\right)  $ in the topology of $H_{0}^{1}\left(  0,1\right)  $ and then in
$C([0,1])$ as well \cite{ARV}. Therefore, for the given $\varepsilon_{0}$ $:=
$ $\min_{x\in\lbrack x_{0},x_{1}]}\left(  y\left(  x\right)  -v_{1}^{+}\left(
x\right)  \right)  $ there exists $T=T\left(  \varepsilon_{0}\right)  $ such
that%
\[
dist_{C([0,1])}(G(t,u_{0}),\mathcal{A})\leq\frac{\varepsilon_{0}}{2}\text{ if
}t\geq T.
\]
In particular, this implies that%
\[
u_{s}\left(  t,x\right)  \leq v_{1}^{+}\left(  x\right)  +\frac{\varepsilon
_{0}}{2}\text{ if }t-s\geq T\text{, }x\in\lbrack0,1]\text{.}%
\]
We choose $s_{0}$ such that $t_{0}-s_{0}\geq T$. Then $y=\gamma\left(
t_{0}\right)  \leq u_{s_{0}}\left(  t_{0}\right)  $ gives%
\[
y\left(  x\right)  \leq u_{s_{0}}\left(  t_{0},x\right)  \leq v_{1}^{+}\left(
x\right)  +\frac{\varepsilon_{0}}{2}<y\left(  x\right)  \text{ for }%
x\in\lbrack x_{0},x_{1}]\text{,}%
\]
which is a contradiction.

Finally, by the same argument in the proof of Lemma \ref{upper}
$\widetilde{\gamma}(t)=-\gamma\left(  t\right)  $ is a bounded complete orbit
if $\gamma\left(  t\right)  $ is a bounded complete orbit. Hence,
$\widetilde{\gamma}\left(  t\right)  \leq v_{1}^{+}$ implies that%
\[
\gamma\left(  t\right)  \geq-v_{1}^{+}=v_{1}^{-}\text{ for all }t\in
\mathbb{R}.
\]
This proves that $y\geq v_{1}^{-}$ for all $y\in\mathcal{A}\left(  t\right)
,\ t\in\mathbb{R}.$
\end{proof}

\bigskip

\begin{theorem}
\label{Structure2}There exist bounded complete trajectories $\gamma_{\ast
}\left(  t\right)  ,\gamma^{\ast}(t)$ such that:

\begin{enumerate}
\item $v_{1}^{-}\leq\gamma_{\ast}\left(  t\right)  \leq\gamma^{\ast}\left(
t\right)  \leq$\ $v_{1}^{+}$ for $t\in\mathbb{R}.$

\item $\mathcal{A}(t)\subset I_{\gamma_{\ast}}^{\gamma^{\ast}}\left(
t\right)  $ for all $t\in\mathbb{R}.$

\item $\gamma_{\ast}$ ($\gamma^{\ast}$) is minimal (maximal) in the sense
given in Theorem \ref{Structure}.

\item If, in addition, $b,\omega\in W_{loc}^{1,2}(\mathbb{R})$, then
$\gamma_{\ast}$ ($\gamma^{\ast}$) is globally asymptotically stable from below (above).
\end{enumerate}
\end{theorem}

\begin{proof}
We shall check the conditions of Theorem \ref{Structure}. By Lemma \ref{upper}
and Theorem \ref{Intervalo}, property (\ref{AttrInterval}) is satisfied with
$\overline{u}\left(  t\right)  \equiv v_{1}^{+},\ \underline{u}\left(
t\right)  \equiv v_{1}^{-}$. From the pullback attraction of $\{\mathcal{A}%
(t)\}_{t\in\mathbb{R}}$, (\ref{AttrSubUpper}) obviously holds. Thus, the
conditions of Theorem \ref{Structure} hold and the results of the first three
items follow. Finally, we observe that the fact that the solution
corresponding to the initial datum $\gamma_{\ast}(s)$ ($\gamma^{\ast}(s)$) is
unique for any $s\in\mathbb{R}$ if $b,\omega\in W_{loc}^{1,2}(\mathbb{R})$ was
proved in \cite[Lemma 22]{Valero2021}. Hence, item 4 follows again from
Theorem \ref{Structure}.
\end{proof}

\begin{lemma}
\label{GammaStar}$\gamma_{\ast}\left(  t\right)  =-\gamma^{\ast}\left(
t\right)  $ for all $t\in\mathbb{R}.$
\end{lemma}

\begin{proof}
As $\varphi\in\mathcal{R}_{\tau}$ implies that $-\varphi\in\mathcal{R}_{\tau}%
$, we obtain for any bounded complete orbit $\gamma\left(
\text{\textperiodcentered}\right)  $ that%
\[
-\gamma\left(  t\right)  \leq\gamma^{\ast}\left(  t\right)  \text{ for all
}t\in\mathbb{R}.
\]
Hence,%
\[
\gamma_{\ast}\left(  t\right)  \geq-\gamma^{\ast}\left(  t\right)  \ \text{for
all }t\in\mathbb{R}.
\]
Then, from the minimality of $\gamma_{\ast},$ we conclude that $\gamma_{\ast
}\left(  t\right)  =-\gamma^{\ast}\left(  t\right)  $ for all $t\in
\mathbb{R}.$
\end{proof}

\begin{corollary}
$\gamma^{\ast}\left(  t\right)  \geq0,\ \gamma_{\ast}\left(  t\right)  \leq0$
for all $t\in\mathbb{R}.$
\end{corollary}

\bigskip

We define the subset of continuous functions on $[0,1]$ satisfying zero
Dirichlet boundary conditions and being strictly positive on $\left(
0,1\right)  $:%
\[
\Phi(0,1)=\{u\in C([0,1]):u(x)>0\text{, }\forall x\in\left(  0,1\right)
,\ u\left(  0\right)  =u\left(  1\right)  =0\}.
\]
Further, we will prove that $\gamma^{\ast}\left(  t\right)  \in\Phi(0,1)$ for
all $t\in\mathbb{R}$. For this aim, let us also consider the autonomous
problem%
\begin{equation}
\left\{
\begin{array}
[c]{l}%
\dfrac{\partial u}{\partial t}-\dfrac{\partial^{2}u}{\partial x^{2}}\in
b_{0}H_{0}(u)+\omega_{0}u,\text{ on }\left(  0,1\right)  \times(\tau
,\infty),\\
u(t,0)=u(t,1)=0,\\
u(\tau,x)=u_{\tau}(x),
\end{array}
\right.  \label{Aut2}%
\end{equation}
that is, when $b\left(  t\right)  \equiv b_{0},\ \omega\left(  t\right)
\equiv\omega_{0}$. We recall the following comparison principle between
solutions to (\ref{Incl}) and (\ref{Aut2}) \cite[Theorem 3]{CLV2020}.

\begin{lemma}
\label{Comparison3}For any initial datum $u_{\tau}\geq0$ there exists a
non-negative strong solution $\overline{u}\left(  \text{\textperiodcentered
}\right)  $ to (\ref{Incl})\ with $\overline{u}\left(  \tau\right)  =u_{\tau}$
such that%
\begin{equation}
\overline{u}\left(  t\right)  \geq u_{b_{0,}\omega_{0}}\left(  t\right)
\text{, }\forall t\geq\tau, \label{IneqNonautAut3}%
\end{equation}
where $u_{b_{0,}\omega_{0}}\left(  \text{\textperiodcentered}\right)  $ is an
arbitrary strong non-negative solution to (\ref{Aut2})\ with $u_{b_{0,}%
\omega_{0}}\left(  \tau\right)  $ $=$ $u_{\tau}.$
\end{lemma}

Denote now by $v_{1,b_{1,}\omega_{1}}^{+},\ v_{1,b_{0,}\omega_{0}}^{+}$ the
positive stationary points of problems (\ref{Aut}) and (\ref{Aut2}), respectively.

\begin{lemma}
\label{AcotGammastar}$v_{1,b_{0,}\omega_{0}}^{+}\leq\gamma^{\ast}\left(
t\right)  \leq v_{1,b_{1,}\omega_{1}}^{+}$ for all $t\in\mathbb{R}$.
\end{lemma}

\begin{proof}
We only need to prove the first inequality, as the second one was established
in Theorem \ref{Structure2}.

Let $s_{n}\rightarrow-\infty$. By Lemma \ref{Comparison3} there exists a
sequence of non-negative strong solutions $u_{n}\left(
\text{\textperiodcentered}\right)  \in\mathcal{R}_{s_{n}}$ to problem
(\ref{Incl})\ such that $u_{n}\left(  s_{n}\right)  =v_{1,b_{0,}\omega_{0}%
}^{+}$ and $u_{n}\left(  r\right)  \geq v_{1,b_{0,}\omega_{0}}^{+}$ for all
$r\geq s_{n}$. We take an arbitrary $t\in\mathbb{R}.$ As $dist\left(
u_{n}\left(  t\right)  ,\mathcal{A}(t)\right)  \rightarrow0$ as $n\rightarrow
\infty,$ passing to a subsequence we have that
\[
u_{n}\left(  t\right)  \rightarrow y\in\mathcal{A}(t),
\]%
\[
v_{1,b_{0,}\omega_{0}}^{+}\leq y\leq\gamma^{\ast}\left(  t\right)  .
\]

\end{proof}

\bigskip

We have given an alternative proof of the results proved already in
\cite{CLV2020}. Additionally, it was proved in \cite[Corollary 7]{CLV2020}
that $\gamma^{\ast}\left(  \text{\textperiodcentered}\right)  $ is the unique
bounded complete orbit such that $\gamma^{\ast}\left(  t\right)  \in\Phi(0,1)$
for all $t\in\mathbb{R}$, that is, it is the unique non-degenerate bounded
complete trajectory.

This result was improved in \cite{Valero2021}, where it was proved that
$\gamma^{\ast}\left(  \text{\textperiodcentered}\right)  $ is the unique
bounded complete trajectory which is non-degenerate at $-\infty$ under the
additional assumption $b,\omega\in W_{loc}^{1,2}(\mathbb{R})$. We recall that
a complete orbit $\gamma\left(  \text{\textperiodcentered}\right)  $ is
non-degenerate at $-\infty$ if for some $t_{0}$ we have $\gamma(t)\in
\Phi(0,1)$ for all $t\leq t_{0}$. Moreover, as remarked before, in this case
the solution starting at any $\gamma^{\ast}\left(  s\right)  $ is unique.

\subsection{An application of Theorem \ref{abstract_upper_semicont XX}: the
asymptotically autonomous case}

\label{sec32}

Here we will consider problem (\ref{Incl}) with $b:\mathbb{R}\rightarrow
\mathbb{R}^{+},$ $\omega:\mathbb{R}\rightarrow\mathbb{R}^{+}$ continuous
functions such that $0<\beta_{0}\leq b\left(  t\right)  \leq b_{1}$ and
$0\leq\gamma_{0}\leq\omega\left(  t\right)  \leq\omega_{1}<\pi^{2}.$ We assume
that $b(t)\rightarrow b_{0}\geq\beta_{0}$ and $\omega\left(  t\right)
\rightarrow\omega_{0}\geq\gamma_{0}$ as $t\rightarrow+\infty.$ In this
especial case problem (\ref{Incl}) is asymptotically autonomous. We will show
the upper semicontinuity of the pullback attractor with respect to the
attractor of the autonomous limit inclusion, that is, that $\lim
_{t\rightarrow+\infty}\mbox{\rm
dist}(\mathcal{A}(t),\mathcal{A})=0,$ where $\mathcal{A}$ is the global
attractor of the autonomous problem:%
\begin{equation}
\left\{
\begin{array}
[c]{l}%
\dfrac{\partial u}{\partial t}-\dfrac{\partial^{2}u}{\partial x^{2}}\in
b_{0}H_{0}(u)+\omega_{0}u,\text{ on }\left(  0,\infty)\times(0,1\right)  ,\\
u(t,0)=u(t,1)=0,\\
u(0,x)=u_{0}(x).
\end{array}
\right.  \label{Aut3}%
\end{equation}

\begin{theorem}
\label{technical lemma}Let $\tau_{n}\nearrow+\infty.$ If $u_{\tau_{n}}$ $\in$
$\mathcal{A}(\tau_{n})$ and $u_{\tau_{n}}\rightarrow u_{0}$ in $L^{2}(0,1)$ as
$n$ $\rightarrow$ $+\infty,$ then for each family of strong solutions
$u^{n}\left(  \text{\textperiodcentered}\right)  $ of problem (\ref{Incl})
with $u^{n}\left(  \tau_{n}\right)  =u_{\tau_{n}}$ there exists a strong
solution $v\left(  \text{\textperiodcentered}\right)  $ of problem
(\ref{Aut3}) such that, up to a subsequence, $v^{n}(t):=u^{n}(t+\tau_{n})$
$\rightarrow$ $v(t)$ in $L^{2}(0,1),$ as $n$ $\rightarrow$ $+\infty,$
uniformly on compact sets of $[0,+\infty).$
\end{theorem}

\begin{proof}
Let $u^{n}\left(  \text{\textperiodcentered}\right)  $ be strong solutions of
(\ref{Incl}) with $u^{n}(\tau_{n})=u_{\tau_{n}}\in\mathcal{A}(\tau_{n}).$ Then
there exist $f^{n}\in L_{loc}^{2}(\tau_{n},+\infty;L^{2}(0,1))$ such that
$f^{n}(t,x)\in H_{0}(u^{n}(t,x))$ a.e. in $(\tau,+\infty)\times\left(
0,1\right)  $ and%
\begin{equation}
\dfrac{du^{n}}{dt}(t)-Au^{n}(t)=b(t)f^{n}(t)+\omega(t)u^{n}(t),\quad\text{a.e
in}\;(\tau_{n},+\infty), \label{P}%
\end{equation}

From the invariance of the pullback attractor we have that
\[
u^{n}(t+\tau_{n})\in U(t+\tau_{n},\tau_{n},\mathcal{A}(\tau_{n}))=\mathcal{A}%
(t+\tau_{n}).
\]
Since $\cup_{\tau\in\mathbb{R}}\mathcal{A}(\tau)$ is bounded in $H_{0}%
^{1}\left(  0,1\right)  $, there exists a constant $C>0$ such that%
\begin{equation}
\left\Vert u^{n}(t+\tau_{n})\right\Vert _{H_{0}^{1}}\leq C,\text{\ }\forall
t\geq0\text{ and }n. \label{Ineq0}%
\end{equation}

Let us fix an arbitrary $T>0$. From the definition of $H_{0}$ it follows that
$\left\vert f^{n}\left(  t,x\right)  \right\vert \leq1$, so in particular the
sequence $\{g^{n}\left(  \text{\textperiodcentered}\right)  \}$, defined by
$g^{n}($\textperiodcentered$)=f^{n}($\textperiodcentered$+\tau_{n}),$ is
bounded in $L^{\infty}(0,T;L^{2}(0,1))$. Hence, up to a subsequence,
$g^{n}\rightarrow g$ weakly star in $L^{\infty}(0,T;L^{2}(0,1)) $ and weakly
in $L^{2}(0,T;L^{2}(0,1))$ for some function $g\left(
\text{\textperiodcentered}\right)  .$

Multiplying (\ref{P}) by $-\frac{\partial^{2}u}{\partial x^{2}}$ and using
Young's inequality we obtain
\begin{equation}
\frac{d}{dt}\Vert u^{n}\Vert_{H_{0}^{1}}^{2}+2\Vert\frac{\partial^{2}u^{n}%
}{\partial x^{2}}\left(  s\right)  \Vert_{L^{2}}^{2}\leq2b_{1}^{2}\Vert
f^{n}(u^{n}\left(  s\right)  )\Vert_{L^{2}}^{2}+2\omega_{1}\left\Vert
u^{n}\left(  s\right)  \right\Vert _{L^{2}}^{2}+\Vert\frac{\partial^{2}u^{n}%
}{\partial x^{2}}\left(  s\right)  \Vert_{L^{2}}^{2}. \label{Ineq1}%
\end{equation}
Hence,%
\begin{equation}
\int_{\tau_{n}}^{\tau_{n}+T}\Vert\frac{\partial^{2}u^{n}}{\partial x^{2}%
}\left(  s\right)  \Vert_{L^{2}}^{2}ds\leq K_{1}\text{, }\forall n,
\label{Ineq2}%
\end{equation}
for some $K_{1}>0$. Therefore, from equality (\ref{P}) we also get that%
\begin{equation}
\int_{\tau_{n}}^{\tau_{n}+T}\Vert\frac{du^{n}}{ds}\left(  s\right)
\Vert_{L^{2}}^{2}ds\leq K_{2}\text{, }\forall n, \label{Ineq3}%
\end{equation}
for some $K_{2}>0.$

Let $v^{n}\left(  \text{\textperiodcentered}\right)  =u^{n}\left(
\text{\textperiodcentered}+\tau_{n}\right)  $, $b^{n}\left(
\text{\textperiodcentered}\right)  =b\left(  \text{\textperiodcentered}%
+\tau_{n}\right)  ,\ \omega^{n}\left(  \text{\textperiodcentered}\right)
=\omega\left(  \text{\textperiodcentered}+\tau_{n}\right)  $. Then for each
$n$ the function $v^{n}\left(  \text{\textperiodcentered}\right)  $ is the
unique strong solution of the problem%
\[
\left\{
\begin{array}
[c]{l}%
\dfrac{\partial v}{\partial t}-\dfrac{\partial^{2}v}{\partial x^{2}}%
=b^{n}(t)g^{n}(t)+\omega^{n}(t)v(t),\text{ on }(0,T)\times\left(  0,1\right)
,\\
v(t,0)=v(t,1)=0,\\
v(0,x)=u_{\tau_{n}}(x).
\end{array}
\right.
\]
Then (\ref{Ineq0}), (\ref{Ineq2}), (\ref{Ineq3}) imply the existence of a
function $v\left(  \text{\textperiodcentered}\right)  $ and subsequence of
$\{v^{n}\}$ such that%
\[
v^{n}\rightarrow v\text{ weakly star in }L^{\infty}(0,T;H_{0}^{1}(0,1)),
\]%
\[
v^{n}\rightarrow v\text{ weakly in }L^{2}(0,T;D(A)),\text{ }%
\]%
\[
\frac{dv^{n}}{dt}\rightarrow\frac{dv}{dt}\text{ weakly in }L^{2}%
(0,T;L^{2}(0,1)).
\]
The functions $v^{n}:[0,T]\rightarrow L^{2}(0,1)$ are then equicontinuous. As
$v^{n}\left(  t\right)  $ is relatively compact in $L^{2}(0,1)$ for each
$t\in\lbrack0,T]$, the Ascoli-Arzel\`{a} theorem gives that%
\[
v^{n}\rightarrow v\text{ in }C([0,T],L^{2}(0,1)).
\]
Also, the convergences $b\left(  t\right)  \rightarrow b_{0},\ \omega\left(
t\right)  \rightarrow\omega_{0},$ as $t\rightarrow+\infty,$ imply that
\[
b^{n}\rightarrow b_{0},\ \omega^{n}\rightarrow\omega_{0}\text{ in }C([0,T]).
\]
Thus, passing to the limit in (\ref{P}) we obtain that%
\[
\frac{dv}{dt}-Av=b_{0}g+\omega_{0}v\text{ in }L^{2}(0,T;L^{2}(0,1)),
\]
so%
\[
\dfrac{dv}{dt}(t)-Av(t)=b_{0}g(t)+\omega_{0}v(t)\text{ in }L^{2}(0,1)\text{
for a.a. }t\in\left(  0,T\right)  .
\]

Since $v:[0,T]\rightarrow L^{2}(0,1)$ is absolutely continuous, it remains to
check that $g\left(  t,x\right)  \in H_{0}(v(t,x))$ for a.a. $\left(
t,x\right)  .$ Observe first that for a.a. $\left(  t,x\right)  \in\left(
0,T\right)  \times\left(  0,1\right)  $ there exists $N(t,x)$ such that
$g^{n}(t,x)\in H_{0}(v(t,x))$ for any $n\geq N.$ Indeed, let $A\subset
\lbrack0,T]\times\lbrack0,1]$ be a set of zero measure such that
$v^{n}(t,x)\rightarrow v(t,x)$ for any $\left(  t,x\right)  \in A^{c}$. If
$(t_{0},x_{0})\in A^{C}$ satisfies $v\left(  t_{0},x_{0}\right)  =0$, then the
result is evident from $g^{n}(t_{0},x_{0})\in\lbrack-1,1]=H_{0}(v(t_{0}%
,x_{0}))$ for all $n.$ If $(t_{0},x_{0})\in A^{C}$ is such that $v(t_{0}%
,x_{0})>0$, then $v^{n}(t_{0},x_{0})\rightarrow v(t_{0},x_{0})$ implies the
existence of $N(t_{0},x_{0})$ for which $v^{n}(t_{0},x_{0})>0$ for $n\geq N$,
so $g^{n}(t_{0},x_{0})=1=H_{0}(v(t_{0},x_{0}))$. The same argument is valid
for $v(t_{0},x_{0})<0$. Now, since $g^{n}\rightarrow g$ weakly in
$L^{1}(0,T;L^{2}(0,1))$ and the set $H_{0}(v(t,x))$ is convex, it follows from
Lemma 32 in \cite{Valero2021} that $g\left(  t,x\right)  \in H_{0}(v(t,x))$
for a.a. $\left(  t,x\right)  $.

We have proved that $v\left(  \text{\textperiodcentered}\right)  $ is a strong
solution of problem (\ref{Aut3}) in the interval $[0,T]$. If we consider the
sequence $T,$ $2T,\ 3T$ and apply a diagonal argument we obtain the desired
strong solution $v\left(  \text{\textperiodcentered}\right)  $ defined in
$[0,+\infty)$ and that%
\[
u^{n}(t+\tau_{n})=v^{n}(t)\rightarrow v(t)
\]
uniformly on compact sets of $[0,+\infty).$
\end{proof}

\bigskip

Using Theorem~\ref{technical lemma} and
Theorem~\ref{abstract_upper_semicont XX} we obtain

\begin{theorem}
\label{uppersemicontinuity} Let $\{\mathcal{A}(t)\}_{t\in\mathbb{R}}$ be the
pullback attractor of problem (\ref{Incl}) and $\mathcal{A}$ the global
attractor of the autonomous problem (\ref{Aut3}). We have
\[
\lim_{t\rightarrow+\infty}\mbox{\rm
dist}(\mathcal{A}(t),\mathcal{A})=0.
\]

\end{theorem}

Finally, we will prove that the complete orbits $\gamma^{\ast}\left(
\text{\textperiodcentered}\right)  ,\ \gamma_{\ast}\left(
\text{\textperiodcentered}\right)  $ converge to the positive and negative
equilibria, respectively, of the autonomous problem (\ref{Aut3}).

\begin{theorem}
$\gamma^{\ast}\left(  t\right)  \rightarrow v_{1}^{+}$ as $t\rightarrow
+\infty$, $\gamma_{\ast}\left(  t\right)  \rightarrow v_{1}^{-}$ as
$t\rightarrow+\infty.$
\end{theorem}

\begin{proof}
Take a sequence $\tau_{n}\nearrow+\infty$. We note that $\gamma^{\ast}\left(
\tau_{n}\right)  \in\mathcal{A}(\tau_{n})$ and that passing to a subsequence
$\gamma^{\ast}\left(  \tau_{n}\right)  \rightarrow\varphi_{0}$. By Theorem
\ref{technical lemma} there exists $\psi^{0}\in\mathcal{R}_{0} $ such that the
sequence $v_{n}^{0}($\textperiodcentered$):=\gamma^{\ast}( $%
\textperiodcentered$+\tau_{n})$ satisfies, passing to a subsequence, that
$v_{n}^{0}(t)\rightarrow$ $\psi^{0}(t)$ in $L^{2}(0,1),$ as $n$ $\rightarrow$
$+\infty,$ for all $t\geq0$. Now, we consider the sequence $v_{n}^{-1}\left(
\text{\textperiodcentered}\right)  =\gamma^{\ast}($\textperiodcentered
$+\tau_{n}-1)$. Hence, there exist $v^{1}\left(  \text{\textperiodcentered
}\right)  \in\mathcal{R}_{0}$ and a subsequence (of the previous subsequence)
such that $v_{n}^{-1}\left(  t\right)  \rightarrow v^{1}\left(  t\right)  $
for all $t\geq0$. We define the function $\psi^{1}:[-1,+\infty)\rightarrow
L^{2}\left(  0,1\right)  $ given by $\psi^{1}\left(  t\right)  =v^{1}\left(
t+1\right)  $. Then from the previous convergence it is clear that $\psi
^{1}\left(  t\right)  =\psi^{0}\left(  t\right)  $ for all $t\geq0$. If we
continue in this way, we construct a sequence of functions $\psi
^{k}:[-k,+\infty)\rightarrow L^{2}\left(  0,1\right)  $ satisfying that
$\psi^{k}(t)=\psi^{k-1}(t)$, for all $t\geq-k+1,$ and that $v^{k}\left(
\text{\textperiodcentered}\right)  =\psi^{k}\left(  \text{\textperiodcentered
}-k\right)  \in\mathcal{R}_{0}$. Defining $\psi:\mathbb{R}\rightarrow
L^{2}\left(  0,1\right)  $ to be equal to the common value of the functions
$\psi^{k}$ at any time $t$, we obtain a bounded complete orbit of
$\mathcal{R}_{0}$. Also, by a diagonal argument there is a subsequence such
that $\gamma^{\ast}(t+\tau_{n})\rightarrow\psi\left(  t\right)  $ for all $t.$

Let us prove that $\psi\left(  t\right)  \equiv v_{1}^{+}$. By Lemma
\ref{AcotGammastar} we know that $v_{1,\beta_{0},\gamma_{0}}^{+}\leq
\gamma^{\ast}\left(  t\right)  $, where $v_{1,\beta_{0},\gamma_{0}}^{+}$ is
the positive equilibrium of problem (\ref{Aut2}) with $b_{0}=\beta_{0}$,
$\omega_{0}=\gamma_{0}$. Thus, $v_{1,\beta_{0},\gamma_{0}}^{+}\leq\psi\left(
t\right)  $ and then $\psi\left(  t\right)  \in\Phi\left(  0,1\right)  $ for
all $t$, so it is a bounded complete non-degenerate orbit. As remarked in the
previous section, $\gamma^{\ast}\left(  \text{\textperiodcentered}\right)  $
is the unique bounded complete non-degenerate orbit for the nonautonomous
problem (\ref{Incl}), and the same is true for the positive equilibrium
$v_{1}^{+}$ in the autonomous counterpart. Hence, $\psi\left(  t\right)
\equiv v_{1}^{+}$.

We need to show that $\gamma^{\ast}\left(  t\right)  \rightarrow v_{1}^{+}$ as
$t\rightarrow+\infty$. If not, there would exist $\varepsilon>0$ and a
sequence $\tau_{n}\nearrow+\infty$ such that $\rho\left(  \gamma^{\ast}\left(
\tau_{n}\right)  ,v_{1}^{+}\right)  >\varepsilon$ for all $n$. By the previous
argument, there exists a subsequence $\{\gamma^{\ast}\left(  \tau_{n_{k}%
}\right)  \}$ such that $\gamma^{\ast}\left(  \tau_{n_{k}}\right)  \rightarrow
v_{1}^{+}$, which is a contradiction.

Finally, $\gamma_{\ast}\left(  t\right)  \rightarrow v_{1}^{-}$ follows from
the fact that $\gamma_{\ast}\left(  t\right)  =-\gamma^{\ast}\left(  t\right)
\rightarrow-v_{1}^{+}=v_{1}^{-}.$
\end{proof}

\bigskip

\section{Final remarks}

We guess that it is also possible to establish the main result for a general
universe $\mathcal{D}$ instead of only bounded sets (see
\cite{LangaSuarez2002} for the single-valued situation).

\bigskip

\textbf{Acknowledgments}

This work has been partially supported by the Spanish Ministry of Science,
Innovation and Universities, project PGC2018-096540-B-I00, by the Spanish
Ministry of Science and Innovation, project PID2019-108654GB-I00, by the Junta
de Andaluc\'{\i}a and FEDER, project P18-FR-4509, by the Generalitat
Valenciana, project PROMETEO/2021/063, and by the Fondo Europeo de Desarrollo
Regional (FEDER) and Consejer\'{\i}a de Econom\'{\i}a, Conocimiento, Empresas
y Universidad de la Junta de Andaluc\'{\i}a, Programa Operativo FEDER
2014-2020
$\backslash$%
\# US-1254251 and
$\backslash$%
\# P20-00592.

Also, J. Simsen and M.S. Simsen have been partially supported by FAPEMIG,
projects APQ-01601-21 and APQ-00675-21.


\bigskip

\begin{thebibliography}{99}                                                                                               %


\bibitem {aragao}E.R. Arag\~{a}o-Costa, T. Caraballo, A.N. Carvalho, and J.A.
Langa, \textit{Stability of gradient semigroups under perturbation},
Nonlinearity, Vol. 24 (2011), 2099-2117.

\bibitem {aragao2}E.R. Arag\~{a}o-Costa, T. Caraballo, A.N. Carvalho, and J.A.
Langa, \textit{Continuity of Lyapunov functions and of energy level for a
generalized gradient system}, Topol. Methods Nonlinear Anal., Vol. 39 (2012), 57-82.

\bibitem {ArnoldChu}L. Arnold, and I. Chueshov, \textit{Order-preserving
random dynamical systems: equilibria, attractors, applications}, Dynam.
Stability Systems, Vol. 13\textbf{\ }(1998), 265-280.

\bibitem {ARV}J.M. Arrieta, A. Rodr\'{\i}guez-Bernal, and J. Valero,
\textit{Dynamics of a reaction-diffusion equation with a discontinuous
nonlinearity}, Internat. J. Bifur. Chaos Appl. Sci. Engrg., Vol.
16\textbf{\ (}2006), 2695-2984.

\bibitem {BV}A.V. Babin, and M.I. Vishik, "Attractors of Evolution Equations,"
North Holland, Amsterdam, 1992.

\bibitem {BCCL}M.C. Bortolan, T. Caraballo, A.N. Carvalho, and J.A. Langa,
\textit{Skew-product semiflows and Morse decomposition}, J. Differential
Equations, Vol. 255 (2013), 2436-2462.

\bibitem {BCLBook}M.C. Bortolan, A.N. Carvalho, and J.A. Langa, "Attractors
under autonomous and non-autonomous perturbations," volume 246 of Mathematical
Surveys and Monographs, American Mathematical Society, 2020.

\bibitem {BCLR}M.C. Bortolan, A.N. Carvalho, J.A. Langa, and G. Raugel,
\textit{Non-autonomous perturbations of Morse-Smale semigroups: stability of
the phase diagram}, J. Dynamics Differential Equations, in press, https://doi.org/10.1007/s10884-021-10066-6.

\bibitem {BrCarVal}R.C.D.S. Broche, A.N. Carvalho, and J. Valero, \textit{A
non-autonomous scalar one-dimensional dissipative parabolic problem: the
description of the dynamics}, Nonlinearity, Vol. 32 (2019), 4912-4941.

\bibitem {CLMV2003}T.Caraballo, J.A.Langa, V.S. Melnik, and J.Valero,
\textit{Pullback attractors of nonautonomous and stochastic multivalued
dynamical systems}, Set-Valued Anal., Vol. 11\textbf{\ }(2003), 153-201.

\bibitem {CLV2005}T. Caraballo, J.A. Langa, and J. Valero, \textit{Asymptotic
behaviour of monotone multi-valued dynamical systems}, Dyn. Syst., Vol.
20\textbf{\ }(2005), 301-321.

\bibitem {CLV2016}T. Caraballo, J.A. Langa, and J. Valero, \textit{Structure
of the pullback attractor for a non-autonomous scalar differential inclusion},
Discrete Contin. Dyn. Syst. Ser. S, Vol. 9\textbf{ }(2016), 979-994.

\bibitem {CLV2020}T. Caraballo, J.A. Langa, and J. Valero, \textit{Extremal
bounded complete trajectories for nonautonomous reaction-diffusion equations
with discontinuous forcing term}, Rev. Mat. Complut., Vol. 33 (2020), 583-617.

\bibitem {CL09}A.N. Carvalho, and J.A. Langa, \textit{An extension of the
concept of gradient semigroups which is stable under perturbation}, J.
Differential Equations, Vol. 246 (2009), 2646-2668.

\bibitem {CLRBook}A.N. Carvalho, J.A. Langa, and J.C. Robinson, "Attractors
for infinite-dimensional non-autonomous dynamical systems," volume 182 of
Applied Mathematical Sciences, Springer, New York, 2013.

\bibitem {CLRS}A.N. Carvalho, J.A. Langa, J.C. Robinson, and A. Su\'{a}rez,
\textit{Characterization of non-autonomous attractors of a perturbed gradient
system}, J. Differential Equations, Vol. 236 (2007), 570-603.

\bibitem {CLR12}A.N. Carvalho, J.A. Langa, and J.C. Robinson,
\textit{Structure and bifurcation of pullback attractors in a non-autonomous
Chafee-Infante equation}, Proc. Amer. Math. Soc., Vol. 140 (2012), 2357-2373.

\bibitem {CI}N.~Chafee, and E.F. Infante, \textit{A bifurcation problem for a
nonlinear partial differential equation of parabolic type}, Applicable Anal.,
Vol. 4 (1974/75), 17-37.

\bibitem {CV}V.V. Chepyzhov, and M.I. Vishik, "Attractors for equations of
mathematical physics," Colloquium Publications 49, American Mathematical
Society, Providence, R.I., 2002.

\bibitem {conley}C. Conley, "Isolated invariant sets and the Morse index,"
CBMS Regional Conference Series in Mathematics Vol. 38, American Mathematical
Society, Providence, R.I., 1978.

\bibitem {Chueshov}I. Chueshov, \textit{Order-preserving skew-product flows
and nonautonomous parabolic systems}, Acta Appl. Math., Vol. 65 (2001), 185-205.

\bibitem {CV16}H.B. da Costa, and J. Valero, \textit{Morse decompositions and
Lyapunov functions for dynamically gradient multivalued semiflows}, Nonlinear
Dynam., Vol. 84 (2016), 19-34.

\bibitem {Ha}J.K. Hale, "Asymptotic behavior of dissipative systems,"
Mathematical Surveys and Monographs V. 25, American Mathematical Society,
R.I., Providence, 1988.

\bibitem {HR}J.K. Hale, and G. Raugel, \textit{Lower semi-continuity of
attractors of gradient systems and applications}, Ann. Mat. Pur. Appl., Vol.
154 (1989), 281-326.

\bibitem {HR92}J.K. Hale, and G. Raugel, \textit{Convergence in dynamically
gradient systems with applications to PDE}, Z. Angew. Math. Phys., Vol. 43
(1992b), 63-124.

\bibitem {HMO}J.K. Hale, L.T. Magalh\~{a}es, and W.M. Oliva, "An introduction
to infinite-dimensional dynamical systems - geometric theory," Applied
Mathematical Sciences Vol.\ 47, Springer-Verlag, 1984.

\bibitem {He}D. Henry, "Geometric theory of semilinear parabolic equations,"
Lecture Notes in Mathematics Vol. 840, Springer-Verlag, Berlin, 1981.

\bibitem {Henry85}D. Henry, \textit{Some infinite-dimensional Morse-Smale
systems defined by parabolic partial differential equations}, J. Differential
Equations, Vol. 59 (1985), 165-205.

\bibitem {HS}M.W. Hirsch, and H. Smith, "Monotone Dynamical Systems, Handbook
of Differential Equations: Ordinary Differential Equations," Volume 2, Chapter
4, North-Holland, 2006, pp.239-357.

\bibitem {KapKasVal11}O.V. Kapustyan, P.O Kasyanov, and J. Valero,
\textit{Pullback attractors for a class of extremal solutions of the 3D
Navier-Stokes system}, J. Math. Anal. Appl., Vol. 373\textbf{\ }(2011), 535-547.

\bibitem {KapKasVal14}O.V. Kapustyan, P.O Kasyanov, and J. Valero,
\textit{Structure and regularity of the global attractor of a
reacction-diffusion equation with non-smooth nonlinear term}, Discrete Contin.
Dyn. Syst., Vol. 34 (2014), 4155-4182.

\bibitem {La}O.A. Ladyzhenskaya, "Attractors for semigroups and evolution
equations," Cambridge University Press, Cambridge, 1991.

\bibitem {LangaSuarez2002}J.A. Langa, and A. Suarez, \textit{Pullback
permanence for non-autonomous partial differential equations}, Electron. J.
Differential Equations Vol. 2002 (2002) (No. 72), pp.1-20.

\bibitem {Rob}J.C. Robinson, "Infinite-dimensional dynamical systems,"
Cambridge University Press, Cambridge, 2001.

\bibitem {RRV}J.C. Robinson, A. Rodr\'{\i}guez-Bernal, and A. Vidal-L\'{o}pez,
\textit{Pullback attractors and extremal complete trajectories for
non-autonomous reaction-diffusion problems}, J. Differential Equations, Vol.
238 (2007), 289-337.

\bibitem {RBV}A. Rodr\'{\i}guez-Bernal, and A. Vidal-L\'{o}pez,
\textit{Extremal equilibria for reaction-diffusion equationsin bounded domains
and applications}, J. Differential Equations, Vol. 244 (2008), 2983-3030.

\bibitem {SY}G.R. Sell, and Y. You, "Dynamics of evolutionary equations,"
Applied Mathematical Sciences Vol. 143, Springer-Verlag, New York, 2002.

\bibitem {SS}J. Simsen, and M.S. Simsen, \textit{On asymptotically autonomous
dynamics for multivalued evolution problems}, Discrete Contin. Dyn. Syst. Ser.
B, Vol. 24 (2019), 3557-3567.

\bibitem {SV}J. Simsen, and J. Valero, \textit{Characterization of pullback
attractors for multivalued nonautonomous dynamical systems}, In: V.A.
Sadovnichiy, M.Z. Zgurovsky Eds., Advances in dynamical systems and control,
Studies in Systems, Decision and Control 69, Springer, Cham, 2016, pp. 179-195.

\bibitem {Te}R. Temam, "Infinite-dimensional dynamical systems in mechanics
and physics," Springer-Verlag, Berlin, 1967.

\bibitem {Valero2012}J. Valero, A weak comparison principle for
reaction-diffusion systems, J. Funct. Spaces Appl., Vol. 2012 (2012), 30pp.

\bibitem {Valero2021}J. Valero, \textit{Characterization of the attractor for
nonautonomous reaction-diffusion equations with discontinuous nonlinearity},
J. Differential Equations, Vol. 275 (2021), 270-308.

\bibitem {Vi92}M.I. Vishik, "Asymptotic behaviour of solutions of evolutionary
equations," Cambridge University Press, Cambridge, 1992.
\end{thebibliography}
\end{document}